\def\rest{\mathord{\restriction}}

\renewcommand{\phi}{\varphi}

\newcommand{\sat}{\models}
\newcommand{\su}{\subseteq}
\renewcommand{\a}{\alpha}
\renewcommand{\b}{\beta}
\newcommand{\g}{\gamma}

\renewcommand{\d}{\delta}
\renewcommand{\l}{\lambda}
\renewcommand{\k}{\kappa}

\newcommand{\om}{\omega}

\newcommand{\lng}{\langle}
\newcommand{\rng}{\rangle}
\newcommand{\ov}{\overline}

\newcommand{\name}[1]{\underset{\sim}{#1}}

\newcommand{\acc}{{\operatorname {acc}}}

\newcommand{\dom}{{\operatorname {dom}}}
\newcommand{\cf}{{\operatorname {cf}}}

\newcommand{\otp}{{\operatorname {otp}}}
\newcommand{\ran}{{\operatorname  {ran}}}

\newcommand{\On}{{\text {On}}}

\newcommand{\Cal}[1]{{\mathcal{#1}}}
\newcommand{\force}{\Vdash}
\newcommand{\imply}{\Rightarrow}
\newcommand{\comp}{\text{comp}\,}
\newcommand{\col}{\text{Col}\,}

\newcommand{\N}{\mathbb N}

\documentclass{amsart}
 \usepackage{amstext}
 \usepackage{amssymb}

\newtheorem{lemma}{Lemma}[section]
\newtheorem{preltheorem}[lemma]{Theorem}
\newtheorem{theorem}[lemma]{Theorem}

 \newtheorem{definition}[lemma]{Definition}  
  \newtheorem{corollary}[lemma]{Corollary}

 \newtheorem{claim}[lemma]{Claim}  
 \newtheorem{fact}[lemma]{Fact}

\usepackage{graphics}
\renewcommand{\int}{\operatorname{\rm int}}

\begin{document}

\title{Fallen Cardinals}

\author{Menachem Kojman}
\address{Department of Mathematics and Computer Science\\
Ben Gurion University of the Negev\\
Beer Sheva, Israel}

\email{kojman@math.bgu.ac.il}

\author{Saharon Shelah}
\address{Institute of Mathematics\\
The Hebrew University of Jerusalem\\
Jerusalem 91904, Israel} 

\email{shelah@math.huji.ac.il} 

\thanks{The second author was partially supported by the Israeli 
Foundation for Basic Science.  Number 720 in list of publications}

\date{December, 1999}


\begin{abstract} 
We prove that for every singular cardinal $\mu$ of cofinality $\om$,
the complete Boolean algebra $\comp\Cal P_\mu(\mu)$ contains
 as a complete
subalgebra an isomorphic copy of the collapse algebra
$\text{Comp}\,\text{Col}(\om_1,\mu^{\aleph_0})$.  Consequently, adding a
generic filter to the quotient algebra $\Cal P_\mu(\mu)=\Cal
P(\mu)/[\mu]^{<\mu}$ collapses
$\mu^{\aleph_0}$ to $\aleph_1$.  Another corollary is that the Baire 
number of the space $U(\mu)$ of all uniform ultrafilters over $\mu$ is 
equal to $\om_2$. The corollaries affirm two  conjectures by Balcar
and Simon. 

The proof uses pcf theory. 
\end{abstract}
\maketitle

\section{Introduction}
\subsection{Forcing and distributivity of complete Boolean algebras}
Every separative poset $P$ which may be used as a forcing notion, is
embedded as a dense subset of a (unique) complete Boolean algebra,
called the \emph{completion of $P$} and denoted by $\comp P$.  The
properties of the forcing extension $V^P$ of the universe $V$ of set
theory, which is obtained by forcing with $P$, are tightly related to
the Boolean-algebraic properties of $\comp P$, in particular to the
\emph{distributivity} properties of $\comp P$.  The least cardinality
of a new set in $V^P$, for example, is equal to the
\emph{distributivity number} of $\comp P$, denoted $\frak h(\comp P)$,
which should really be called the ``non-distributivity number'', since
it is the least cardinality of a product of sums which violates
distributivity (see \cite{JechHB} for more information).  Finer
properties of non-distributivity determine which cardinals of $V$ are
preserved and which are collapsed in the extension $V^P$.  The
non-distributivity property which is important for our context is the
following:

\begin{definition}
A complete Boolean algebra $B$ is $(\k,\cdot,\l)$-nowhere distributive
iff $B$ contains partitions of unity $P_\a$  for $\a<\k$ (namely,
$\sum P_\a=1$ and $p\wedge q=0$ for $p\not=q$ in $P_\a$)
 so that for every $b\in B-\{0\}$ there exists $\a<\k$ so
that $b\wedge p\not= 0$ for $\ge \l$ members $q\in P_\a$.
\end{definition}

Clearly, if $\l_1<\l_2$ and $B$ is $(\k,\cdot,\l_2)$ nowhere
distributive, it is also $(\k,\cdot,\l_1)$-nowhere distributive.  The
systematic study of distributivity in Boolean algebras was pursued by
the Czech school of set theory ever since the discovery of Forcing in
1963.

 \subsection{The quotient algebra $\Cal P_\mu(\mu)$} In 1972 Balcar
 and Vop\v enka began the study non-distributivity in quotient algebras
 $\Cal P_\k(\k)$ for infinite cardinals $\k$.  For every infinite
 cardinal $\k$, the algebra $\Cal P_\k(\k)$ is obtained as the
 quotient of the power set algebra $\Cal P(\k)$ over the ideal
 $[\k]^{<\k}$ of all subsets of $\k$ whose cardinality is strictly
 smaller than $\k$.  It was first shown that for every singular
 cardinal $\mu$ of countable cofinality the distributivity number of
 $\Cal P_\mu(\mu)$ is $\om_1$ and that for every cardinal $\k$ of
 uncountable cofinality the distributivity number of $\Cal P_\k(\k)$
 is $\om$ \cite{BVop}.  The distributivity number of $\Cal P_\om(\om)$
 was discussed separately in \cite{BPS}.  The exact nature of
 distributivity in various $\Cal P_\k(\k)$ was addressed in a series
 of papers \cite{BF,BSHB,Balcar}, usually under additional set
 theoretic assumptions.  The optimal ZFC non-distributivity properties
 of $\Cal P_\k(\k)$ were obtained in \cite{Balcar}, from which we
quote:

\begin{theorem}[Balcar and Simon]\label{nwd}
\begin{enumerate}
\item For every singular $\mu$ of countable cofinality $\Cal
P_\mu(\mu)$ is $(\om_1,\cdot,\mu^{\aleph_0})$-nowhere distributive
\item For every singular $\k$ of uncountable cofinality $\Cal
P_\k(\k)$ is $(\om,\cdot,\k^+)$-nowhere distributive.
\end{enumerate}
\end{theorem}

\subsection{The collapse algebra $\col (\om_1,\mu^{\aleph_0})$} It was
in \cite{BF} that it was first shown that under certain set theoretic
assumptions $\Cal P_\k(\k)$ (and some other factor algebras of $\Cal
P(\k)$) have completions which are isomorphic to suitable
\emph{collapse algebras}.  Let us introduce collapse algebras.  For
cardinals $\k<\l$, $\k$ regular, the poset $\col(\k,l)$ is the natural
$\k$-complete poset for introducing a function $\varphi$ from $\k$
onto $\l$, namely for ``collapsing'' $\l$ to $\k$.

\begin{equation}
    \begin{split}
     	\col(\k,\l)=\{f:\text{ for some $\a< \k$,  } f \text{ is a functions, }\\
 	\dom f=\a \text{ and } \ran f\su \l\}
	\end{split}
 \end{equation}
 
 The completion $\comp \col(\k,\l)$ is the collapse algebra for
 $(\k,\l)$.  The cardinality of $\col(\k,\l)$ is clearly $\l^{<\k}$,
 and therefore, since $\col (\k,\l)$ is dense in its completion, the
 density $\pi(\comp\col(\k,\l))$ is equal to $\l^{<k}$.  For each
 $\a<\k$ let $P_\a$ be a maximal antichain in $\col(\k,\l)$ composed
 of conditions which decide $\varphi\rest\a$.  The set
 $\{P_\a:\a<\k\}$ (which is, really, a \emph{name} for $\phi$) is also
 a witness to the fact that $\comp\col(\k,\l)$ is
 \emph{$(\k,\cdot,\l)$-nowhere distributive}.  The following
 characterization of $\comp\col(\om_1,\mu^{\aleph_0})$ for a singular
 $\mu$ of countable cofinality is a particular instance of a general
 characterization theorem for collapse algebras (\cite{BSHB}, 1.15):

\begin{theorem}\label{charac} Let $B$ be a complete
$(\om_1,\cdot,\mu)$ nowhere distributive Boolean algebra containing an
$\aleph_1$-closed dense subset.  If $\pi(B)=\mu^{\aleph_0}$, then $B$
is isomorphic to $\comp\col(\om_1,\mu^{\aleph_0})$.
\end{theorem}

Thus $\comp\col(\om_1,\mu^{\aleph_0})$ is characterized by the
existence of an $\aleph_1$-closed dense subset, density
$\mu^{\aleph_0}$ and $(\om_1,\cdot,\mu)$-nowhere distributivity (which
is a weaker condition than $(\om_1,\cdot,\mu^{\aleph_0})$-nowhere
distributivity which is actually satisfied).

 \subsection{The problem}  Let $\mu$ be a singular cardinal of
countable cofinality.  Two of the properties which characterize
$\comp\col(\om_1,\mu^{\aleph_0})$ hold also in the quotient algebra
$\Cal P_\mu(\mu)$: $\aleph_1$-completeness (easily)  and
$(\om_1,\cdot,\l)$-nowhere distributivity (by Theorem \ref{nwd}). 
Could it be  true that $\comp\,\Cal P_\mu(\mu)$ and
$\comp\;\text{Col}\;(\om_1,\mu^{\aleph_0})$ are isomorphic?  An old
independence result of Baumgartner's rules that out.  Baumgartner
forced an almost disjoint family in $\Cal P(\mu)$ of size
$2^\mu>\mu^{\aleph_0}$, showing thus that it is consistent with ZFC
that the cellularity, hence density, of $\Cal P_\mu(\mu)$ strictly
exceeds $\mu^{\aleph_0}$ (\cite{Baumgartner}, 6.1).  In Baumgartner's
model $\text{Comp}\,\Cal P_\mu(\mu)$ cannot be isomorphic to
$\text{Comp}\,\text{Col}\,(\om_1,\mu^{\aleph_0})$, whose density is
exactly $\mu^{\aleph_0}$.

However, if one assumes that $\mu^{\aleph_0}=2^\mu$, it follows
trivially that the density of $\Cal P_\mu(\mu)$ is $\mu^{\aleph_0}$,
and hence, by the aforementioned characterization of $\comp
\col(\om_1,\mu^{\aleph_0})$, it is isomorphic to $\comp\Cal
P_\mu(\mu)$.   In particular,   denoting by $\vdash$  provability in ZFC, we
have
\cite{Balcar}:

\begin{equation}\label{conj1}
 \mu^{\aleph_0}=2^\mu\vdash \mu^{\aleph_0} \text{
collapses to } \aleph_1 \text{ in } V^{\Cal P_\mu(\mu)}
\end{equation}

What, then, is the precise relation between $\Cal P_\mu(\mu)$ and
$\col\;(\om_1,\mu^{\aleph_0})$?  Most importantly, does forcing with
$\Cal P_\mu(\mu)$ always collapse $\mu^{\aleph_0}$ to $\aleph_1$?

Balcar and Simon conjectured in \cite{Balcar1} that the answer is
``yes", namely, that the cardinal arithmetic assumption
$\mu^{\aleph_0}=2^{\mu}$ could be removed from (\ref{conj1}).  In the
same paper they advance towards an affirmative solution of their
conjecture by proving in ZFC that forcing with $\Cal P_\mu(\mu)$
collapse the continuum $2^{\aleph_0}$ to $\om_1$.  Since for
$\mu<2^{\aleph_0}$ it holds trivially that $\mu^{\aleph_0}=
2^{\aleph_0}$, that proves their conjecture for all countably cofinal
singular cardinals $\mu$ which are below the continuum:

\begin{equation}\label{cont}
\mu<2^{\aleph_0}\vdash 
\mu^{\aleph_0} \text{ collapses to } \aleph_1 \text{ in } V^{\Cal
P_\mu(\mu)}
\end{equation}

Finally, there was the problem of computing the \emph{Baire number} of
the space $U(\mu)$ of all uniform ultrafilters over $\mu$.  An
ultrafilter $F$ over $\mu$ is uniform if it does not contain a set of
cardinality $<\mu$.  With the usual topology, in which the basic open
sets are $\hat p=\{F\in U(\mu):p\in F\}$ for $p\in [\mu]^\mu$, the
space $U(\mu)$ is a compact Hausdorff space and is therefore not
coverable by $\om_1$ nowhere-dense sets.  The Baire number of a space
with no isolated points is the least number of nowhere-dense sets
needed to cover the space.  In \cite{Balcar1} it was proved that the
Baire number of $U(\mu)$ is $\om_2$ under any of the following
assumptions: $(i)$ $2^{\aleph_0}>\aleph_1$, $(ii)$
$2^\mu=\mu^{\aleph_0}$ or $(iii)$ $2^{\om_1}=\om_2$.  It was
conjectured that the Baire number of $U(\mu)$ could be shown to equal
$\om_2$ in ZFC alone.

\subsection{The solution} 
 The main result in the present paper determines the precise relation
 between $\text{Comp}\,\Cal P_\mu(\mu)$ and
 $\text{Comp}\,\text{Coll}\;(\om_1,\mu^{\aleph_0})$. 
 The collapse algebra  is isomorphic to a
complete
 subalgebra of the quotient algebra (Theorem \ref{main} below):
\begin{equation}\label{emb}
     \vdash \comp\col(\om_1,\mu^{\aleph_0}) \lessdot \comp\Cal P_\mu(\mu)
    \end{equation}

This implies that the universe $V^{\Cal P_\mu(\mu)}$ contains
$V^{\text{Coll}\;(\om_1,\mu^{\aleph_0})}$ as a subuniverse. 
Therefore,
\begin{equation}\label{final}
  \vdash  \mu^{\aleph_0} \text{ collapses to } \om_1\text{ in }V^{\Cal
P_\mu(\mu)}
\end{equation}

which proves the conjecture. An easy corollary of (\ref{emb}) is that the Baire
number of $U(\mu)$ is equal to $\om_2$ (Theorem \ref{baire} below).

Balcar and Simon stated in \cite{Balcar1} another ZFC conjecture concerning
singular cardinals  of uncountable cofinality.
The authors will present a
solution of that conjecture in a sequel paper. 

\subsection{History} B.~Balcar presented this conjecture to the
authors during a meeting in Hattingen, Germany, in June of 1999. 
Shelah then proved, using the Erd\H os-Rado theorem, that:

\begin{equation}\label{shelah}
   \mu > 2^{\aleph_0}\vdash  \mu^{\aleph_0}\text{ collapses 
    to } 2^{\aleph_0}\text{  in } V^{\Cal P_\mu(\mu)}
    \end{equation}
    
This affirmed the ZFC conjecture, since (\ref{cont}) and
(\ref{shelah}) together give (\ref{final}) ($2^{\aleph_0}=\mu$ is of course
impossible by K\" onig's Lemma).  In August of 1999 Kojman found a  ZFC proof
of (4) by replacing Shelah's use of the Erd\H os-Rado theorem (which requires
cardinal arithmetic assumptions) by a use of a pcf theorem.  This
proof is presented below. 

\subsection{Description of the proof} 
Let $P=\lng [\mu]^\mu,\le\rng$ where, for $p_1,p_2\in P$, $p_1\le
p_2\iff |p_1-p_2|<\mu$.  For every $D\su \Cal P_\mu(\mu)$, $D$ is
dense in $\Cal P_\mu(\mu)$ if and only if $\bigcup D$ is dense in $P$
and $D$ is a filter in $\Cal P_\mu(\mu)$ if and only if $\bigcup D$ is
a filter in $P$.  Therefore, $G\su \Cal P_\mu (\mu)$ is a generic
filter over $\Cal P_\mu(\mu)$ if and only if $\bigcup G$ is a generic
filter over $P$.  Hence $V^{\Cal P_\mu(\mu)}=V^P$.  For convenience,
we work with $P$ rather than with $\Cal P_\mu(\mu)$.  Let $\l$ denote
$\mu^{\aleph_0}$.

The main point in finding a complete copy of $\comp\col(\om_1,\l)$
inside $\comp\Cal P_\mu(\mu)$ is to overcome the large cellularity
that $\text{Comp}\,P$ may possess, e.g. in Baumgartner's model.  This
is achieved by forcing only with $Q\su P$, which contains all
\emph{closed} conditions of $P$.  It is not hard to verify that
$\text{Comp}\,Q$ is isomorphic to a complete subalgebra of
$\text{Comp}\,P$.  Then it is shown that
$\text{Comp}\,Q\cong\text{Comp}\,\text{Col}(\om_1,\l)$.  To that end
one needs to prove that $\pi(Q)=\l$.  This fact is achieved by an old
trick: club guessing.  Once density is out of the way, it remains to
establish $(\om_1,\cdot,\mu)$-nowhere distributivity of $\comp Q$, to
facilitate the use of  Theorem
\ref{charac} above.  Here another pcf tool is used: the Trichotomy
Theorem.

\subsection{Notation and preliminaries} Our notation is mostly
standard. One exception is that when  the relations $f<_U g$, $f\le_U g$
for ordinal functions $f,g$ where $U$ an ultrafilter over $\om$
is extended to \emph{partial} functions. We recall that if
$P$ and $Q$ are posets and for some $P$-name $\name G$ it holds that
$\force_P$ ``$\name G$ is a generic filter over $Q$" and for every
$q\in Q$ there exists $p\in P$ such that $p\force q\in \name G$, then
$\text{Comp}\, Q$ is isomorphic to a complete subalgebra of the
$\text{Comp}\,P$ via the embedding $b\mapsto \sum\{p\in P: p\force_P
``b\in \name G''\}$.

 The following  two theorems from pcf
theory will be used:

\begin{preltheorem}[Club Guessing] \label{clubguess} If $\k^+<\l$ and
$\k,\l$ are regular cardinals, then there exists a sequence $\ov
C=\lng c_\d:\d<\l \wedge \cf\d=\k\}$ so that:

\begin{enumerate}
\item For every $\d<\l$ with $\cf\d=\k$, $c_\d$ is closed and
unbounded in $\d$ and $\otp\, c_\d=\k$.  \item For every club $E$ of
$\l$ there exists $\d\in S^\l_\k$ so that $c_\d\su E$.
\end{enumerate}
\end{preltheorem}

\begin{preltheorem} \label{trichotomy}(The Trichotomy)
Suppose $A$ is an infinite set, $I$ an ideal over $A$ and $\l>|A|^+$ a
regular cardinal. If $\ov f= \lng f_\a:\a<\l\rng$ is a
$<_I$-increasing sequence of ordinal functions on $A$, then one of the
following conditions holds:
\begin{itemize}
\item (Good) $\ov f$ has an exact upper bound $f$ with $\cf f(a)>|A|$ for
all $a\in A$;

\item (Bad)  there are sets $S(a)$ for $a\in A$ satisfying $|S(a)|\le
|A|$ and an ultrafilter $D$ over $A$ extending the dual of $I$ so that
for all $\a<\l$ there exists $h_\a\in\prod S(a)$ and $\b<\l$ such that
$f_\a <_D h_\a <_D f_\b$.

\item (Ugly) there is a function $g:A\to \On$ such that letting
$t_\a=\{a\in A: f_\a(a)>g(a)\}$, the sequence $\ov t=\lng
t_\a:\a<\l\rng$ does not stabilize modulo $I$.

\end{itemize}
\end{preltheorem}

A proof of the club guessing Theorem can be found in \cite{cardarith},
III,\S 1, \cite{ABC} or the appendix to \cite{409}.  The Trichotomy
Theorem is Lemma 3.1 in \cite{cardarith}, and a shorter proof of it is
available in the appendix to \cite{EUB}.

\subsection{Acknowledgments} The authors wish to thank Bohuslav Balcar
for presenting this problem to them during the ESF meeting in
Hattingen, June 1999, and take the opportunity to thank R\" udiger G\"
obel and Simone Pabst for their wonderful work in organizing this
meeting. The result in this paper was presented by  Kojman in the last of
 five lectures on pcf theory which were delivered at the Winter School on
General Topology held in the Czech republic in January 2000. His   thanks  for
the hospitality and for the opportunity to present this material to an
interested audience are  extended here to the organizers.

\section{The proof}

Throughout this Section let $\mu$ be a fixed singular cardinal
of countable cofinality, and let $\mu_n$ be a fixed  strictly increasing
sequence of regular cardinals with $\sum \mu_n=\mu$.  Let $P=\{p\su
\mu:|p|=\mu\}$ and for $p_1,p_2\in P$ let $p_1\le p_2$  iff
$|p_1-p_2|<\mu$.  As pointed out in the introduction, forcing with $P$
is equivalent to adding a generic filter to $\Cal P_\mu(\mu)$.

Denote $\l=\mu^{\aleph_0}=|\prod \mu_n|$. 

\begin{theorem}\label{main} 
    $\comp\col(\om_1,\l) \lessdot \comp\Cal P_\mu(\mu)$.
\end{theorem}

\begin{proof}

Let $Q=\{q\in P: q \text{ is closed} \}$.  For a condition $p\in P$,
let $\acc\, p$ be the set of all accumulation points of $p$.  Clearly,
$\acc\, p\in Q$

\begin{lemma}\label{Q}
Let $G\su P$ be a generic filter.  Then $G_1=\{q\in Q: (\exists p\in
G)(\acc\,p\le q\}$ is generic in $Q$.
\end{lemma}

\begin{proof} 
If $p_1,p_2\in G$ then there is some $p_3\in G$ so that $p_3\su
p_1\cap p_2$.  So $\acc\, p_3\le \acc\, p_1,\acc\, p_2$ and $\acc\,
p_3\in\{\acc\, p:p\in G\}$.  Thus $G_1$ is closed under finite
intersections.  Clearly, $G_1$ is upwards closed.  Thus $G_1$ is a
filter.

Suppose that $D\su Q$ is dense and downwards closed.  Let $p\in P$ be
arbitrary, and consider $q=\acc\, p$.  Let $q_1\le q$ be chosen in
$D$.

For $\a\in q_1$ define $\b_\a=\min(p-(\a+1))$, and let
$p_1=\{\b_\a:\a\in q_1\}$.  $q_2=\acc\, p_1\su \acc\, q_1\le q_1$ so
there is some $p_1\le p$ with $\acc\, p_1\in D$.
\end{proof}

By this Lemma it follows that 
\begin{equation}\label{Qcomp}
\comp Q\lessdot \comp P
\end{equation}

We aim now to show that \begin{equation}\label{Q=Col}
\comp Q\cong \comp\col(\om_1,\l)
\end{equation} 

First, we shall see that 
$\pi(Q)=\l$. 

Let $q\in Q$ be arbitrary. Let $a(p)=\{n:p\cap
[\mu_n,\mu_{n+1})\not=\emptyset\}$ and let $\{m_n:n<\om\}$ be the
increasing enumeration of $a(p)$.  

\begin{definition}
A condition $q\in Q$ is \emph{normal} if it satisfies 
\begin{equation}\label{normal}
                   \otp\,[q\cap [\mu_{m_n},\mu_{m_n+1})]=\mu_n+1		   
\end{equation}
\end{definition}

\begin{lemma}
    The set of normal conditions is dense in $Q$.
\end{lemma}

\begin{proof}
Given a condition $q$, let $m_n$ be the least so that $q\cap
[\mu_{m_n},\mu_{m_n+1})|\ge \mu_n$ and choose $c_n\su q\cap
[\mu_{m_n},\mu_{m_n+1})$ of order type $\mu_n + 1$.  Let $q'=\bigcup_n
c_n$.  Thus $q'\le q$, and is of the form (\ref{normal}) above.
\end{proof}

Let $M$ be a fixed  elementary submodel of $\lng H(\Omega),\in\rng$ for
a sufficiently large regular cardinal $\Omega$ so that $\mu\su M$, 
 $[M]^{\aleph_0}\su M$ and the cardinality of $M$ is $\l$.  Let
$Q^M=Q\cap M$. Clearly, $|Q^M|=\l$.

\begin{lemma}\label{density}
$Q^M$ is dense in $Q$.
\end{lemma}

\begin{proof} Let $q$ be a condition in $Q$ and assume, without loss
of generality, that it is normal.  Let $c_n=q\cap
[\mu_{m_n},\mu_{m_n+1})$.

\begin{claim}For every $n$, there exists a closed subset of $c_{n+2}$
 of order type  $\mu_n + 1$ which belongs to $M$.
\end{claim}

\begin{proof}
Let $\g=\sup c_{{n+2}}$.  In $M$, fix an increasing and continuous
function $f:\mu_{n+2}\to \g$ with $\sup \ran f=\gamma$.  Let
$E=\{i<\mu_{{n+2}}:f(i)\in c_{n+2}\}$.  Thus $E\su \mu_{{n+2}}$ is a
club in $\mu_{{n+2}}$.

The club $E$ itself may not belong to $M$ (because $c_{n+2}$ may not
belong to $M$). But since $\mu_n^+<\mu_{n+2}$ and both (regular)
cardinals belong to $M$, $M$ contains some club guessing sequence
$\lng c_\d:\d<\mu_{{n+2}} \wedge \cf\d=\mu_n\rng$ by the club guessing Theorem
\ref{clubguess} above. Thus there is some
$\d<\mu_{n+2}$ so that $c_\d\su E$. Clearly, $c_\d\in M$. Since $f\in
M$, also $\ran(f\rest c_\d)\cup \{\sup \ran(f\rest c_\d)\} \su E$
belongs to $M$, and is a closed subset of $c_{{n+2}}$ of order type
$\mu_n+1$.
\end{proof}

Using the claim, choose, for every $n$, a closed set $b_n$ so that
$b_n\su c_{n+2}$, $\otp\, b_n=\mu_n+1$ and $b_n\in M$.  Since $M$ is
closed under countable sequences, $q'=\bigcup b_n\in Q^M$, and clearly
$q'\le q$ is a normal condition in $Q$.
\end{proof}

 This has established that $\pi(Q)=\l$. 

We  need the following  simple fact about $Q$ and  $Q^M$: 

\begin{fact}\label{QMcomplete}
 $Q$ is $\aleph_1$-complete and $Q^M$ is $\aleph_1$-complete.
\end{fact}

\begin{proof}
Suppose that $q_0\ge q_1\ge\ldots$ is a decreasing sequence of
conditions in $Q$.  By induction on $n$, let $m_{n}$ be chosen so that
$\otp\, [q_{n}\cap [\mu_{m_n},\mu_{m_n+1})]>\sum_{i<n} |q_n-q_i|^+$,
and choose a closed subset $c_{n+1}$ of $\bigcap _{i\le n} q_i\cap
[\mu_{m_n},\mu_{m_n+1})$ with $\otp\, c_n=\mu_n + 1$.  The condition
$\bigcup c_n$ belongs to $Q$ and $q\le q_n$ for all $n$.  If each
$q_n$ belongs to $M$ then the sequence itself belongs to $M$ because
$M$ is closed under taking $\om$-sequences, and hence some $q$ which
satisfies $q\le q_n$ for all $n$ belongs to $M$, by elementarity. 
(Alternatively, one can do the induction for proving completeness of
$Q$ inside $M$).
\end{proof}

Thus, $\comp Q$ contains an $\aleph_1$-complete dense set of size
$\l$.  To prove (\ref{Q=Col}) from Theorem \ref{charac}  it remains to show
that
$\comp Q$ is
$(\om_1,\cdot,\mu)$-nowhere distributive. For this purpose we 
  inspect the generic cut which $Q$ creates in $\prod \mu_n/U$, where 
  $U$ is  the generic ultrafilter over $\om$
which forcing with $Q$ introduces. 

\begin{fact} \label{generic} Suppose $G\su Q$ is a generic filter. Then
$\{a(q):q\in G\}$ is an ultrafilter over $\om$.
\end{fact}

\begin{proof}
 If $q_1\le q_2$ are normal conditions, then $a(q_1)-a(q_2)$ is
 finite.  Thus $a:\{q\in Q: q \text{ is normal} \}\to \Cal P(\om)$ is
 an order preserving map onto $\lng\Cal P(\om),\su^*\rng$. 
 Furthermore, if $t\le a(q)$, then there is some $q'\le q$ such that
 $a (q')=t$.  Therefore the image of a generic $G\su Q$ under $a$ is
 an ultrafilter over $\om$.
\end{proof}

Given a normal condition $q\in Q$, define the following two functions
on $a(q)$ by letting, for each $n\in a(q)$,

\begin{gather}
\chi^+_q(n)= \sup [q\cap [\mu_n,\mu_{n+1})]\notag\\
\chi^-_q(n)=
\min [q\cap [\mu_n,\mu_{n+1})]\notag
\end{gather}

The set of conditions $q$ for which $\chi^+_q\in \prod \mu_n$ is
clearly dense in $Q$, so we always assume that $\chi^+_q\in \prod
\mu_n$.

Since
$\chi^-_q(n)<\chi^-+q(n)$ for every
$n\in a(q)$, and
$q\force ``a(q)\in U$'', it follows that $q\force ``\chi^-_q <_U \chi^+_q$''.

Let $G\su Q$ be a generic filter. Define

\begin{gather}D_0^+=\{\chi^+_q:q\in G \text{ and $q$ is normal}\}\notag\\
D_0^-=\{\chi^-_q:q\in G \text{ and $q$ is normal}\}\notag
\end{gather}

Now for each normal $q\in Q$, $q\force \chi^+_q\in D_0^+\wedge
\chi^-_q\in D_0^-$''.

Let

\begin{gather}
    D^+=\{f\in \prod \mu_n: (\exists g\in D_0^+)(g\le_U f)\}\notag\\
    D^-=\{f\in \prod \mu_n: (\exists g\in D^-_0)(f\le_U g)\}\notag
\end{gather}

\begin{lemma}\label{avoidance}
Suppose that $f\in \prod \mu_n$ and $q\in Q$ is a normal condition. 
Then there exists a normal condition $q'\le q$ so that \[q'\force
\chi^+_{q'}<_U f \vee f<_U \chi^-_{q'}\]
\end{lemma}

\begin{proof} If for some infinite set $B\su A(q)$,  $m_n\in a(q)\imply
f(m_n)<\chi^+_q(m_n)$ then $q':=\bigcup _{m_n\in B} c_n-(f(m_n)+1)$ is a
normal condition and for all $n\in a(q')$ it holds that
$\chi^-_{q'}(n)>f(n)$. Since $q'\force a(q')\in U$, the second
alternative holds for $q'$.

If $\{n\in a(q): f(n)<\chi_q^+(n)\}$ is finite, let $n_0$ be fixed so
that for every $n>n_0$ it holds that $\sup c_n\le f(n)$ and let, for
$n>n_0$, $b_n\su c_n$ be the initial segment of $c_n$ whose order type
is $\mu_{n-1} +1$.  Now $\bigcup b_n$ is a normal condition and
$q'\force \chi^+_q<_U f$.
\end{proof}

Since $Q$ is $\aleph_1$-complete, no new members are added to $\prod
\mu_n$ after forcing with $Q$.  Therefore, by Lemma \ref{avoidance},
it holds that $D^-$ is a lower half of a Dedekind cut in $\prod
\mu_n/U$ whose upper half is $D^+$; that $D^-$ has no last element and
that $D^+$ has no first element.  By the definition of $D^-$, it is
clear that $D^-_0$ is cofinal in $(D^-,<_U)$.
Furthermore, if $\{f_i:i<\om\}$ is a set of functions, $q\in Q$ and
for all $i<\om$ it holds that $q\force f_i\in D^-$ then by iterated
use of Lemma \ref{avoidance} and $\aleph_1$-completeness there exists
$q'\le q$ so that $q'\force \bigwedge_i f_i<_U \chi^-_{q'}$.  As a
consequence, the cofinality of $D^-$ is uncountable.

We shall need the following strengthening of Lemma \ref{avoidance},
which says that the generic cut $(D^-,D^+)$  is not trapped by
any product of countable sets.

\begin{lemma}\label{notbad}
Suppose $A_n\su [\mu_n,\mu_{n+1})$ is a countable set for each
$n<\om$, and $q\in Q$ is a normal condition. Then there is a
condition $q'\le q$ in $Q$ so that for every $n\in a(q')$ it holds
that $A_n\cap (\chi^-_{q'},\chi^+_{q'})=\emptyset$.
\end{lemma}

\begin{proof}
Let $\varepsilon_n<\om_1$ be the order type of $A_n$ and let
$\lng \a^n_i:i<\varepsilon_n\rng$ be the increasing enumeration of
$A_n$. Partition $[\mu_n,\mu_{n+1})$ to the intervals $[\mu_n,\a^n_0)$,
$\{[\a^n_i,\a^n_{i+1}):i+1<\varepsilon_n\}$ and $[\sup
A_n,\mu_{n+1})$. For every $n>0$, choose an interval $I_n$ in the
partition of $[\mu_n,\mu_{n+1})$ so that $|I_n\cap q|=\mu_n$, and let $c_n\su
(I_n\cap q)$ be closed of order type $\mu_{n-1}+1$. Let
$q'=\bigcup I_n\cap q$. 
\end{proof}

\begin{lemma} The cofinality of $D^-$ is $\om_1$.
\end{lemma}

\begin{proof} We have seen that $\cf
(D^-)>\aleph_0$. Suppose now, to the contrary, that $\k>\aleph_1$ is regular
and that $q\force ``\ov f=\lng f_i:i<\k\rng$ is $<_U$ increasing and cofinal in
$D$''.  The Trichotomy Theorem applies to $\ov f$, but:

The third condition ( ``Ugly'') cannot hold , since $U$ is an
ultrafilter.

The first condition  (``Good") cannot hold, because in $D^+$ there is no first
element.

Let us see now that the second condition (``Bad'') cannot hold either. 
Suppose that $q\force ``\lng A_n:n<\om\rng$, witnesses `Bad' for $\ov
f$".  Then,

\begin{equation} \label{bad}
    q\force (\forall \a<\k)(\exists h\in \prod A_n)( \exists \b<\k)
    [f_\a<_U h<_U f_\b]
\end{equation}

By $\aleph_1$-completeness, we may assume that $\lng
A_n:n<\om\rng$ and each $A_n$ belong to the ground model. By Lemma
\ref{notbad} there is a condition $q'\le q$ so that for all $n\in
a(q)$, 

\begin{equation}\label{eight}
A_n\cap \Bigl(\chi^-_{q'}(n),\chi^+_{q'}(n)\Bigr)=\emptyset
\end{equation}

Since $q$ forces that $\ov f$ is cofinal in $D^-$ and $q'\force
\chi^-_{q'}\in D^-\wedge \chi^+_{q'}\in D^+$, there is some $\a<\k$
and $q''\le q'$ so that

\begin{equation}\label{alpha} q''\force \chi^-_{q'}<_U f_\a<_U
\chi^+_{q'}
\end{equation}

By strengthening $q''$ we may assume that for some $\b<\k$ and $h\in
\prod A_n$,

\begin{equation}\label{beta}
q''\force \chi^-_{q'} <_U f_\a<_U h <_U f_\b <_U \chi^+_{q'}
\end{equation}

So there is some $n$ (in fact, infinitely many) so that
\begin{equation}
\chi^-_{q'}(n) < f_\a(n)<h(n)<f_\b(n)<\chi^+_{q'}(n)
\end{equation}

This is a contradiction to (\ref{eight}), since $h(n)\in A_n$.

Thus, the cofinality of $D^-$ is at least $\om_1$ and no more than $\om_1$; so
it is exactly
$\om_1$.
\end{proof}

Since $\cf (D^-,<_U)=\om_1$ and $D^-_0$ is cofinal in $(D^-,<_U)$, using
$\aleph_1$-completeness of $Q$ it is easy to find a sequence of
conditions
$\lng q(i):i<\om_1\rng\su G$ such that
$i<j<\om_1\imply q_(i)\ge q(j)$ and 
$\lng
\chi^-_{q(i)}:i<\om_1\rng$ is ($<_U$-increasing and) cofinal in
$(D^-,<_U)$. Fix a $Q$-name ${\name q}$ for
such a sequence. Observe that if $q_1,q_2$ are incompatible, then $\force
``\neg(q_1\in \ran \name q \wedge q_2\in \ran \name q)"$, since any two
conditions in $\ran \name q$ are compatible. 

\begin{lemma} \label{range} For every $q\in Q$ there is a set
$\{q''_\a:\a<\l\}$ of pairwise incompatible conditions below $q$, so that
for each $\a<\l$ there is $q'_\a\le q$ and $i(\a)$ so that $q'_\a
\force \name q(i(\a))=q''_\a$, and $\{q'_\a:\a<\l\}$ are pairwise incompatible.
\end{lemma}

\begin{proof}
Let $q\in Q$ be a normal condition, and let $c_n=q\cap
[\mu_{m_n},\mu_{m_n+1})$. $c_n$ is  a closed set of order type $\mu_n$. For
each
$n$ let $b_n$ be the initial segment of $c_{n+1}$ of order type
$(\prod _{i\le n}\mu_i,<_{\text{lx}})$, the lexicographic ordering of all
sequences $(x_0,x_1,\dots,x_n)$ in the product
$\mu_0\times\mu_1\dots\times \mu_n$. Identify each member in $b_n$
with the sequence in $\prod_{i\le n}\mu_i$ it corresponds to via the
order isomorphism, and define a projection $\pi_{m,n}:b_n\to b_m$ for
$m<n$ by mapping a sequence of length $n$ to its initial segment of
length $m$. The inverse limit of this system is the set of all
functions  $g\in \prod b_n$ with the property that for all $m<n$,
$\pi_{m,n}g(n)=g(m)$. Denote this set of functions by $L\su
\prod b_n$.

Choose a set of $\l$ different functions $\lng g_\a:\a<\l\rng\su L\cap
\prod \acc\, b_n$ and for each $\a$ let $g_\a'(n)=g_\a(n)+1$.  Let
$q_\a=\bigcup_{n>0} b_n\cap [g_\a(m_n),g'_\a(m_n))$.  Thus each $q_\a$
is a condition below $q$.  Furthermore, if $\a\not=\b$ then from some
point $n_0$ on, either $\chi^+_{q_\a}(n)<\chi^-_{q_\b}(n)$ or
$\chi^+_{q_\b}(n)<\chi^-_{q_\a}(n)$.  Thus $\{q_\a:\a<\l\}$ is a set
of pairwise incompatible conditions below $q$.

For each $\a<\l$,
\begin{align}
  q_\a & \force ``\chi^-_{q_\a}\in D^- \wedge \chi^+_{q_\a}\in D^+
  "\notag\\
   q_\a & \force ``(\exists i<\om_1)( \forall
   j<\om_1)[i<j\imply\chi^-_{q_\a} <_U \chi^-_{\name q{}(j)}]"\notag
\end{align}

Fix $q'_\a\le q_\a$ so that for some $i(\a)<\om_1$ and $q''_\a\le
q_\a$, $q'_\a\force \name {q}{}(i(\a))=q''_\a$.  For $\a<\b<\l$,  since
$q''_\a \le q'_\a\le q_\a$, 
$q''_\b\le q'_\b\le q_\b$ and
$q_\a,q_\b$ are incompatible,  $q'_\a$ is oncompatible with $q'_\b$ and $q''_\a$ is
incompatible with 
$q''_\b$. 
\end{proof}

Fix, for each $i<\om_1$, a maximal antichain $P_i\su Q$ of
conditions that decides $\name {q(\a)}$.  
\begin{claim}
 For every condition $q\in Q$ there exists   some
$i<\om_1$ so that
$q$ is compatible with
$\ge \mu$  members of $P_i$.\end{claim} 
\begin{proof}
Let $q\in Q$ be an arbitrary
condition.  By Lemma \ref{range} there are $\l$ pairwise incompatible
conditions $\{q''_\a:\a<\l\}$ below
$q$, each of which is forced to be $\name q{}(i(\a))$ for some $i(\a)<\om_1$,
by some extension $q'_\a\le q$, and $\{q'_\a:\a<\l\}$ are
pairwise incompatible extensions of $q$.   Since
$\l>\mu$, there is necessarily some fixed $i<\om_1$ so that
$|\{\a<\l:\a(i)=i\}|>\mu$. 
Since  different $q'_\a,q'_\b$ in this set force
different values for $\name q(i)$, they cannot be compatible with the
same member of $P_i$. Thus $q$ is compatible with $\ge \mu$ members of $P_i$.
\end{proof} 

    The last claim  established
$(\om_1,\cdot,\mu)$-nowhere distributivity of $\comp Q$.  By Theorem
\ref{charac} $\comp\col(\om_1,\l)\cong \comp Q$,   and since $\comp Q\lessdot
\comp
\Cal P_\mu(\mu)$, the proof is complete.

\end{proof}
\begin{corollary}\label{finalcollapse}
 $V^{\Cal P_\mu(\mu)}\sat |\l|=\aleph_1$.
\end{corollary}

\begin{proof}
Since $\comp\col(\om_1,\l)$ is a complete subalgebra of $\comp\Cal
P_\mu(\mu)$, the universe $V^{\col(\om_1,\l)}$ is contained in
$V^{\Cal P_\mu(\mu)}$.  Therefore, there is an onto function
$\varphi:\om_1\to \l$ in $V^{\Cal P_\mu(\mu)}$.  Since $\Cal
P_\mu(\mu)$ is $\om_1$-complete, $\om_1$ is preserved in $V^{\Cal
P_\mu(\mu)}$.  Thus, the cardinality of $\l$ in $V^{\Cal P_\mu(\mu)}$
is $\aleph_1$.
\end{proof} 

\begin{corollary}\label{baire}
For every singular cardinal $\mu$ with $\cf \mu=\aleph_0$ the Baire number of
$U(\mu)$, the space of uniform ultrafillters over $\mu$, is equal to $\om_2$.
\end{corollary}

\begin{proof}
By Theorem \ref{main} there exists a dense subset of $Q$ which is
isomorphic to the dense subset $D=\{f:\exists i<\om_1[f:(i+1)\to\l]\}$
of $\col (\om_1,\l)$, namely, there are conditions $\{q_f:r\in D\}\su
Q$ so that $q_f\le q_g \iff g\su f$.

Let $W_{i,\a}=\{q_f:\dom f=i+1 \wedge f(i)=\a\}$.  Define
$V_{i,\a}=\{u\in U(\mu):(\exists p\in u)(\exists q_f\in W_{i,\a})[\acc
p\su q_f]\}$.  It should be clear that $V_{i,\a}$ is a maximal
antichain in $P$.  Let $p\in P$ be arbitrary.  By density of
$\{q_f:f\in D\}$ there exists some $f\in \col(\om_1,\l)$ with domain
$i+1<\om$ so that $q_f\le \acc p$.  By \ref{Q}, for each $\a<\l$ there
exists some $p_\a\le p$ so that $\acc p_\a\le q_{f\cup\lng
i+1,\a\rng}$.  Thus, for every $p\in P$ there exists some $i<\om_1$ so
that $p$ is compatible in $P$ with one member from each $V_{i,\a}$.

Now clearly $O_\b=\bigcup_{i<\om_1,\a\ge \b}$ is dense open in
$U(\mu)$ for each $\b<\l$.  Also, $\bigcap O_\b:\b<\om_2=\emptyset$. 
Thus $U(\mu)$ is coverable by $\om_2$ nowhere-dense sets.  Since it is
known \cite{Balcar1} that $U(\mu)$ cannot be covered by fewer than
$\om_2$ nowhere-dense sets, its Baire number is equal to $\om_2$.
\end{proof}

\section{Concluding remarks}
We first remark that the part of the proof between Lemma
\ref{QMcomplete} and Lemma \ref{range} can be applied verbatim to $P$
instead of to $Q$ to show that $\comp\Cal P_\mu(\mu)$ is
$(\om_1,\cdot,\mu^{\aleph_0})$-nowhere distributive, and constitutes
thus an alternative ZFC proof of
$(\om_1,\cdot,\mu^{\aleph_0})$-nowhere distributivity of $\Cal
P_{\mu}(\mu)$ from the Trichotomy theorem.

Next, we remark that Corollary \ref{finalcollapse} can be derived directly,
without invoking Theorem \ref{charac}, as follows: fix a 1-1 function
$f:Q^M\to \l$ and apply Lemma \ref{range} to $Q^M$, observing that the
set $\{q''_\a:\a<\l\}$ belongs so $M$.  Now fix a function $h:\l\to
\l$ such that for every $A\in [\l]^\l\cap M$, $\ran (h\rest A)=\l$. 
The function $h\circ \name q$ is a collapsing function by a simple
density argument.

We devote now a few words to the role of pcf theory in this proof and
in several other proofs.  Pcf theory was developed to provide bounds
on powers of strong limit singular cardinals, or, better, on the
covering numbers of singular cardinals.  The most well known discovery
of the theory is that poset $\lng \Cal P
_{\aleph_0}(\aleph_\om),\supseteq\rng$ of countable subsets of
$\aleph_\om$ ordered by reverse inclusion has a dense subset of size
$<\aleph_{\om_4}$.  In other words: the cardinality of this poset may
be arbitrarily large, but its density is bounded.

\relax From the point of view of pcf theory, powers of regular cardinals are
the ``soft" part of cardinal arithmetic, which envelopes the hard
``skeleton" of powers of singulars that pcf theory addresses --- the
revised power set function pp.  To read more about this philosophy the
reader is referred to \cite{cardarith} (especially the analytical
index.  \S14), \cite{skeptics} and \cite{ABC}.

The proof above is yet another example of the same theme: a complete
subalgebra of density $\mu^{\aleph_0}$ is uncovered inside
$\text{Comp}\, \Cal P_\mu(\mu)$, whose own density may be
$2^\mu>\mu^{\aleph_0}$ in case the power function at regular cardinals
assumes large values.  The powers of regular uncountable cardinals may
be ``peeled off'' from $\comp\Cal P_\mu(\mu)$ by the club-guessing
technique to get to the ``skeleton" $\comp\col(\om_1,\mu^{\aleph_0})$.

Pcf methods are used also in other contexts to show that various
structures on the power set of a singular cardinal contain ``skeletons''
of  bounded cardinality. We quote  the example
\cite{Dowker} of a  Dowker subspace of  cardinality 
$\aleph_{\om+1}$ inside M.~E.~Rudin's  Dowker space \cite{rudin}, whose
cardinality is $(\aleph_\om)^{\aleph_0}$.

Pcf techniques were used for studying collapses of cardinals by
Cummings \cite{james} (see also \cite{EUB}).

Lastly, we remark that while the role of closed unbounded subsets of
\emph{regular} cardinals in combinatorial set theory is so central
that one could not imagine uncountable combinatorics without them, the
proof above shows that also closed subsets of a \emph{singular}
cardinal may be sometimes useful.


\begin{thebibliography}{1}

\bibitem{BF}
B. Balcar and F. Fran\v ek.
\newblock Completion of factor algebras of ideal
\newblock {\em{Proc. Amer. Math. Soc.}} 100(2):205--212, 1987

\bibitem{Balcar}
Bohuslav Balcar and Petr Simon.
\newblock On collections of almost disjoint families.
\newblock {\em {Commentationes Mathematicae Universitatis Carolinae}},
29(4):631--646, 1988.


\bibitem{BSHB}
Bohuslav Balcar and Petr Simon.
\newblock Disjoint refinements
\newblock in {\em{Handbook of Boolean Algebra}} Vol. 2.  North-Holland, 1989. 

\bibitem{Balcar1}
Bohuslav Balcar and Petr Simon.
\newblock Baire number of the space of uniform ultrafilters.
\newblock {\em{Israel J. Math.}} 92:263--272, 1995.

\bibitem{BPS}
B. Balcar, J. Pelant and P. Simon.
\newblock The space of ultrafilters on $\N$ covered by nowhere dense sets.
\newblock {\em Fund. Math.}110:11--24, 1980. 

\bibitem{BVop}
B. Balcar and P. Vop\v enka.
\newblock On Systems of almost disjoint sets.
\newblock {\em{Bull. Acad. Polon. Sci., Ser. Sci Math}} 20:421--424, 1972. 



\bibitem{Baumgartner}
James~E. Baumgartner.
\newblock Almost-disjoint sets, the dense set problem and the partition
  calculus.
\newblock {\em Ann. Math. Logic}, 9(4):401--439, 1976.


\bibitem{james}
James Cummings.
\newblock Collapsing successors of singulars.
\newblock {\em{ Proc. Amer. Math. Soc.}}125:2703--2709, 1997.


\bibitem{JechHB}
Thomas Jech.
\newblock Distributive laws
\newblock in {\em{Handbook of Boolean Algebra}} Vol. 2.  North-Holland, 1989. 


\bibitem{ABC}
Menachem Kojman.
\newblock The abc of pcf.
\newblock {\em Circulated notes}, 1995.

\bibitem{EUB}
Menachem Kojman.
\newblock Exact upper bounds and their uses in set theory.
\newblock {\em Ann. Pure Appl. Logic}, 92(3):267--282, 1998.

\bibitem{409}
Menachem Kojman and Saharon Shelah.
\newblock Non-existence of universal orders in many cardinals.
\newblock {\em Journal of Symbolic Logic}, 57:875--891, 1992.


\bibitem{Dowker}
Menachem Kojman and Saharon Shelah.
\newblock A {Z}{F}{C} {D}owker space in $\aleph\sb {\omega+1}$: an application
  of {P}{C}{F} theory to topology.
\newblock {\em Proc. Amer. Math. Soc.}, 126(8):2459--2465, 1998.

\bibitem{CountTrich}
Menachem Kojman and Saharon Shelah.
\newblock The Condition in the Trichotomy 
Theorem is Optimal.
\newblock {\em Archive for Math. Logic.}, in press. 

\bibitem{rudin} M.~E.~Rudin, {\em A normal space $X$ for which
$X\times I$ is not normal}, Fund. Math. 73 (1971) 189--186

\bibitem{skeptics}
Saharon Shelah.
\newblock {Cardinal arithmetic for skeptics}.
\newblock {\em {Bull. Amer. Math. Soc. New Series}}, {\bf
  26}:197--210, 1992.

\bibitem{cardarith}
Saharon Shelah.
\newblock {\em {Cardinal Arithmetic}}, volume~29 of {\em {Oxford Logic
  Guides}}.
\newblock {Oxford University Press}, 1994.


\end{thebibliography}
\end{document}